\documentstyle{article}
\begin{document}
\title{{Singular Solution to Special Lagrangian
 Equations} }
 \author{{Nikolai Nadirashvili\thanks{LATP, CMI, 39, rue F. Joliot-Curie, 13453
Marseille  FRANCE, nicolas@cmi.univ-mrs.fr},\hskip .4 cm Serge
Vl\u adu\c t\thanks{IML, Luminy, case 907, 13288 Marseille Cedex
FRANCE, vladut@iml.univ-mrs.fr} }}

\date{}
\maketitle

\def\n{\hfill\break} \def\al{\alpha} \def\be{\beta} \def\ga{\gamma} \def\Ga{\Gamma}
\def\om{\omega} \def\Om{\Omega} \def\ka{\kappa} \def\lm{\lambda} \def\Lm{\Lambda}
\def\dl{\delta} \def\Dl{\Delta} \def\vph{\varphi} \def\vep{\varepsilon} \def\th{\theta}
\def\Th{\Theta} \def\vth{\vartheta} \def\sg{\sigma} \def\Sg{\Sigma}
\def\bendproof{$\hfill \blacksquare$} \def\wendproof{$\hfill \square$}
\def\holim{\mathop{\rm holim}} \def\span{{\rm span}} \def\mod{{\rm mod}}
\def\rank{{\rm rank}} \def\bsl{{\backslash}}
\def\il{\int\limits} \def\pt{{\partial}} \def\lra{{\longrightarrow}}

 {\em Abstract.} We prove the existence of non-smooth solutions to
Special Lagrangian Equations in the non-convex case.

\section{Introduction}
\bigskip

In this paper we study a  fully nonlinear second-order elliptic
equations of the form (where $h\in {\bf R} $)
$${\bf F}_h(D^2u)= \det(D^2u)-Tr(D^2u) +h\sigma_2(D^2u)-h=0\leqno(1)$$
defined in a smooth-bordered domain of $\Omega \subset {\bf R}^3$,
$\sigma_2(D^2u)=\lambda_1\lambda_2+\lambda_2\lambda_3
+\lambda_1\lambda_3$ being the second symmetric function of the
eigenvalues
 $\lambda_1,\lambda_2,\lambda_3$
of $D^2u$. Here $D^2u$ denotes the Hessian of the function $u$.
This equation  is equivalent to the  Special Lagrangian potential
equation [HL1]:
$$ \noindent SLE_{\theta}: \:\: \quad \quad\quad \quad\quad \quad\quad \quad
Im  \{e^{-i\theta} \det(I + iD^2u)\}=0 \quad \quad\quad
\quad\quad\quad\quad \quad $$
for $ h:=-\tan(\theta)$ which can be
re-written as
 $${\bf F}_{\theta}=\arctan\lambda_1+\arctan \lambda_2+\arctan\lambda_3-\theta=0 .$$

The set  $$\{A\in Sym^2({\bf R}^3): \:{\bf F}_h(A)=0\}\subset Sym^2({\bf R}^3) $$
 has three connected components, $C_i, i=1,2,3$ which correspond to the values
 $ \theta_1=-\arctan(h)-\pi,\theta_2=-\arctan(h), \theta_3=-\arctan(h)+\pi . $

We study  the Dirichlet problem
$$\cases{{\bf F}_{\theta}(D^2u)=0 &in $\Om$\cr
u=\vph &on $\pt\Om\;,$\cr}$$ where  $\Omega \subset {\bf R}^n$ is
a bounded domain with smooth boundary $\partial \Omega$ and $\vph$
is a continuous function on $\pt\Om$.

\smallskip
For $ \theta_1=-\arctan(h)-\pi$ and $\theta_3=-\arctan(h)+\pi  $
the operator  ${\bf F}_{\theta}$ is concave
 or convex, and the Dirichlet problem  in these cases  was treated in [CNS];
  smooth solutions are
 established  there  for smooth boundary data on appropriately convex domains.

The middle branch $C_2, \theta_2=-\arctan(h) $ is never convex
(neither concave), and the  classical solvability
 of the Dirichlet problem remained open.

In the case of uniformly elliptic equations a theory of weak
(viscosity) solutions for the Dirichlet problem gives  the
uniqueness of such solutions, see [CIL], moreover these solutions
lie in $C^{1,\varepsilon}$ by [C],[T1],[T2]. However, the recent
results [NV1],[NV2],[NV3] show that at least in 12 and more
dimensions the viscosity solution of the Dirichlet problem  for a
uniformly elliptic equation can be singular, even in the case when
the operator depends only on eigenvalues of the Hessian.

 One can define viscosity solutions for  non-uniformly elliptic equations
 (such as $SLE_{\theta}$) as well, but in this case the uniqueness  of viscosity
 solution is not known which makes the use of these solutions less convenient.

Recently a new very interesting approach to degenerate elliptic
equations was suggested by Harvey and Lawson [HL2]. They
introduced a new notion of  a weak solution for the Dirichlet
problem for such equations and proved the existence, the
continuity and the uniqueness of these solutions.

 The main purpose of the present note is to show that  the  classical solvability
  for Special Lagrangian Equations {\it does not} hold.

 More precisely, we show the existence for any $\theta \in ]-\pi/2,\pi/2[ $
 of a small ball  $B\subset{\bf R}^3 $ and of an analytic function $\phi$
 on $\partial B$ for which the unique Harvey-Lawson
  solution $u_{\theta}$ of the Dirichlet problem  satisfies :

\smallskip
$(i)\;\;u_{\theta}\in C^{1,1/3}; $

\smallskip
$(ii)\;\;u_{\theta}\notin C^{1,\delta} $ for $\forall\delta>1/3$.

\smallskip
Our construction use  the Legendre transform for solutions of
${\bf F}_{{1\over h }} (D^2u)=0$ which gives solutions of ${\bf
F}_{h}(D^2u)=0 $; in particular, for $h=0$ it transforms solutions
of $\sigma_2(D^2u)=1$ into solutions of $\det(D^2u)=Tr(D^2u).$
This construction could be of interest by itself.

\smallskip
Finally, we think that the following conjecture is quite plausible:

\smallskip\smallskip
{\bf Conjecture.} {\em Any Harvey-Lawson solution  of
$SLE_{\theta}$    on a ball $B$ lies in $C^1(B)$ (if  $\vph$ is
sufficiently smooth).}

\smallskip
  In the case $\theta=0$ these solutions
lie in $C^{0,1}(B)$ by Corollary 1.2 in [T3].

\section{Harvey-Lawson Dirichlet Duality Theory}
In this section we recall the "Dirichlet duality" theory by Harvey
and Lawson [HL], which establishes (under an appropriate  explicit
geometric assumption on the domain $\Omega $) the existence and
uniqueness of continuous solutions of the Dirichlet problem for
fully nonlinear, degenerate elliptic equations
$${\bf F}(D^2u)=0.\leqno(2)$$ Following  the method by Krylov
[Kr] this theory takes a geometric approach to the equation which
eliminates the operator ${\bf F}$ and replaces it with a closed
subset $F$ of the space $Sym^2({\bf R}^n)$ of real symmetric
$n\times n$ matrices, with the property that $\partial F$ is
contained in $\{{\bf F} = 0\}$. We need only the case when $\Omega
$ is a ball when the geometric assumption is automatically true
and thus we do not discuss it below.

The general set-up of the theory  is the following. Let $F$ be  a
given closed subset of the space of real symmetric matrices
$Sym^2({\bf R}^n)$.  The theory formulates and solves the
Dirichlet problem for the equation
$$Hess_x(u)\in \partial F {\hbox{ for all }} x \in \Omega $$
using the functions of "type $F$", i.e.,
which satisfy
$$Hess_x(u)\in  F {\hbox{ for all }} x .$$
A priori these conditions make sense only for $C^2$ functions $u$.
The theory extends the notion to functions which are only upper
semi-continuous.

 A closed subset $ F\subset Sym^2({\bf R}^n)$  is called a
{\em Dirichlet set} if it satisfies the condition
$$F + {\mathcal P} \subset F$$
 where $${\mathcal P} = \{A \in Sym^2({\bf R}^n): A \ge 0\}$$
 is the subset of non-negative
matrices. This condition corresponds to degenerate ellipticity in
modern fully nonlinear theory; it implies that the maximum of two
functions of type $F$ is again of type $F$ which is the key
requirement for solving the Dirichlet problem. Note that
translates, unions (when closed) and intersections of Dirichlet
sets are Dirichlet sets. The {\em Dirichlet dual} set $\tilde F$
is defined as $$ \tilde F:= -(Sym^2({\bf R}^n)\backslash Int
(F)).$$ By Lemma 4.3 in [HL]  this is equivalent to the condition
$$ \tilde F:= \{A\in Sym^2({\bf R}^n):\forall B \in F, A+B \in\tilde{\mathcal P})),$$
 $ \tilde{\mathcal P}$ being the set  of all quadratic forms except those
 that are negative definite.

An upper semi-continuous (USC) function $u$ is called an {\em
subaffine function} if it verifies locally  the condition:

{\em For each affine function $a$, if $u\le a$ on the boundary of
a ball $B$, then $u \le a$ on $B$}.

Note that a $C^2$-function is subaffine if and only if $Hess(u)$
has at least one non-negative eigenvalue at each point.  An USC
 function $u$ is {\em  of type F} if $u + v$ is
subaffine for all $C^2$-functions $v$ of type $\tilde F$. In other
words, $u$  is of type $F$ if for any "test function" $v \in C^2$
of dual type $\tilde F$, the sum $u + v$ satisfies the maximum
principle.   A function $u$ on a domain is said to be {\em
F-Dirichlet} if $u$ is of type $F$ and $-u$ is of type $\tilde F$.
Such a function $u$ is automatically continuous, and at any point
$x$ where $u$ is $C^2$, it satisfies the condition $$Hess_x(u)\in
\partial F {\hbox{ for all }} x \in \Omega. $$

The main result of the theory (in our restricted setting) is
[HL,Theorem 6.2]

\medskip
{\bf THEOREM  (The Dirichlet Problem)}. {\em Let $B\subset {\bf
R}^n$ be a ball, and let $F$ be a Dirichlet set.   Then for each
$\varphi\in C(\partial B )$, there exists a unique $u \in C(B )$
which is an $F$-Dirichlet function on $B$  and equals $\varphi$
on $\partial B$}.

\smallskip
Besides, one has [HL, Remark 4.9.]:

\medskip
{\bf Proposition (Viscosity Solutions)}. {\em In the conditions of
the theorem u is a viscosity solution of $(2)$.}

\medskip
Krylov's idea [K, Theorem 3.2] permits to reconstruct from $F$ a
{\em canonical form} of the operator ${\bf F}$ such that:

\smallskip
1). $\partial F=\{{\bf F} = 0\};$

\smallskip
2). $ F=\{{\bf F} \ge 0\}.$

\smallskip
It is sufficient to  define

\smallskip

${\bf F}(A):=dist(A, \partial F)$ for $A\in F$;

\smallskip
${\bf F}(A):=-dist(A, \partial F)$ for $A\notin F$.

\smallskip
The operator ${\bf F}$ (in its canonical form) is {\em strictly
elliptic} if for any $A\in F$ there  exists $\delta(A)>0$ s.t.
${\bf F}(A +P)\ge \delta(A)\cdot ||P|| $ for all $P\in {{\mathcal
P}}$, and {\em uniformly elliptic} if
 ${\bf F}(A +P)\ge \delta\cdot||P|| $ for all $P\in {{\mathcal P}}$, $A\in F$ and an
 absolute constant $\delta>0$ (note that for ${\bf F}$ in its canonical form
 ${\bf F}(A +P)-{\bf F}(A)\le ||P||$  by definition). Moreover, the
 function ${\bf F}$ is concave iff $F$ is concave,  and is convex
 iff $\tilde F$ is convex.

Below  we will use the Harvey-Lawson  theory only in the case of
Hessian equations, i.e. when ${\bf F}(A)$ depends only on the
eigenvalues $\lambda_1(A)\le\lambda_2(A)\le\ldots\le\lambda_n(A)$
of $A.$ Then the sets $\{{\bf F}=0\}, F,\tilde F$ are stable under
the action of the orthogonal group $O_n({\bf R})$ by conjugation.
Consider the map
$$Sym^2({\bf R}^n)\longrightarrow D_n\subset{\bf R}^n$$
$$\hskip 1.8 cm A \longmapsto (\lambda_1(A),\lambda_2(A),\ldots,\lambda_n(A)) $$
where
$$D_n:=\{(\lambda_1,\lambda_2,\ldots,\lambda_n)
\in{\bf R}^n:\lambda_1\le\lambda_2\le\ldots\le\lambda_n\}.$$ The
images $\{{\bf F_{\lambda}}=0\}, F_{\lambda},\tilde F_{\lambda}$
of $\{{\bf F}=0\}, F,\tilde F$ determine completely their
preimages. The sets
$$\{{\bf F}_{\Lambda}=0\}:=\bigcup_{\sigma \in S_{n} }\{ {\bf F}_{\sigma(\lambda)}=0\}
 \subset {\bf R}^n, $$
 $$ F_{\Lambda}:=\bigcup_{\sigma \in S_{n}}  F_{\sigma(\lambda)}
 \subset {\bf R}^n, $$
 $$\tilde  F_{\Lambda} :=\bigcup_{\sigma \in S_{n}} \tilde  F_{\sigma(\lambda)}
 \subset {\bf R}^n, $$
 where ${\sigma(\lambda)}:=(\lambda_{\sigma(1)},\lambda_{\sigma(2)},
 \ldots,\lambda_{\sigma(n)})$ are $S_n$-invariant subsets in ${\bf
 R}^n$ which determine $\{{\bf F}=0\}, F$ and $\tilde F$ as well.
 Moreover, by [Ba] (cf. Sec. 3 of [CNS]) the set $F$ is    convex
  iff $F_{\Lambda}$ is   convex.

\section{ Some properties of Special Lagrangian Equations}

In this section we give some properties of the Special Lagrangian
Equation
$${\bf F}_{-c}(D^2u)= \det(D^2u)-Tr(D^2u)
-c\sigma_2(D^2u)+c=0.$$ Note first that the the set $\{{\bf
F}_{-c,\Lambda}=0\}$ is a real cubic surface ${\bf S}_{c,\Lambda}$
with three components ("branches") which can be presented as a
graph:
$$\lambda_3={c(1-\lambda_1\lambda_2)+\lambda_1+\lambda_2\over
\lambda_1\lambda_2-1+c(\lambda_1+\lambda_2)}\;. \leqno (3) $$

One easily proves the following by brute force computations:

\medskip
{\bf Lemma 3.1.} {\em

\smallskip
$1).$ The components of  ${\bf S}_{c,\Lambda}$ are given by
$$C_1=\{(\lambda_1,\lambda_2,\lambda_ 3):\lambda_1\lambda_2-1+
c(\lambda_1+\lambda_2)>0,\; \lambda_1>-c,\lambda_2>-c\}, $$
 $$C_2=\{(\lambda_1,\lambda_2,\lambda_ 3):\lambda_1\lambda_2-1+
 c(\lambda_1+\lambda_2)<0 \}, $$
 $$C_3=\{(\lambda_1,\lambda_2,\lambda_ 3):\lambda_1\lambda_2-1+
 c(\lambda_1+\lambda_2)>0,\;\lambda_1<-c,\lambda_2<-c\}; $$ equivalently,
  $$C_1=\{(\lambda_1,\lambda_2,\lambda_ 3):\arctan\lambda_1+\arctan\lambda_2
  +\arctan\lambda_3= \pi+\arctan c\}, $$
$$C_2=\{(\lambda_1,\lambda_2,\lambda_ 3):\arctan\lambda_1+\arctan\lambda_2+
  \arctan\lambda_3=\arctan c
  \}, $$
$$C_3=\{(\lambda_1,\lambda_2,\lambda_ 3):\arctan\lambda_1+\arctan\lambda_2
  +\arctan\lambda_3=-\pi+\arctan c\;
  \}. $$

\smallskip
  $2).$ For any $c\in {\bf R}$, $C_1$ is  convex,  $C_3$ is concave,  $C_2$  is neither.}

\medskip
   {\em Proof.} 1). is straightforward, 2). follows from the
   Hessian of $\lambda_ 3$ in (3):
   $${\pt^2\lambda_3\over\pt \lambda_1^2}={2(\lambda_2+c)(\lambda_2^2+1)
   (1+c^2)\over(\lambda_1\lambda_2-1+c\lambda_1+c\lambda_2)^3},$$
 $${\pt^2\lambda_3\over\pt \lambda_2^2}={2(\lambda_1+c)(\lambda_1^2+1)
   (1+c^2)\over(\lambda_1\lambda_2-1+c\lambda_1+c\lambda_2)^3},$$
$$\det(D^2\lambda_3)={4(c^2+1)^2(c\lambda_1+c\lambda_2+\lambda_1^2
   \lambda_2^2+\lambda_1\lambda_2+\lambda_2^2+\lambda_1^2)\over
   (\lambda_1\lambda_2-1+c\lambda_1+c\lambda_2)^5},$$
which implies e.g. that the point with $\lambda_1=\lambda_2=
-c-1/10$ for $c\ge 0$,  $\lambda_1=\lambda_2= -c+1/10$ for $c\le
0$ is a saddle point on $C_2$.

\medskip
The corresponding Dirichlet sets $  F^i_{c},\; i=1,2,3$
  are given  (via $F^i_{c,\Lambda}$) by
$$F^1_{c,\Lambda}=\{(\lambda_1,\lambda_2,\lambda_ 3):\arctan\lambda_1+\arctan\lambda_2+
\arctan\lambda_3\ge\pi+\arctan c \}, $$
 $$F^2_{c,\Lambda}=\{(\lambda_1,\lambda_2,\lambda_ 3):\arctan\lambda_1+\arctan\lambda_2+
\arctan\lambda_3\ge\arctan c\}, $$
$$F^3_{c,\Lambda}=\{(\lambda_1,\lambda_2,\lambda_ 3):\arctan\lambda_1+\arctan\lambda_2+
\arctan\lambda_3\ge -\pi+\arctan c \}, $$
 and their duals by ([HL], Prop. 10.4.):
$$\tilde F^1_{c,\Lambda}=\{(\lambda_1,\lambda_2,\lambda_ 3):\arctan\lambda_1
  +\arctan\lambda_2+\arctan\lambda_3\ge-\pi-\arctan c\}, $$
$$ \tilde F^2_{c,\Lambda}=\{(\lambda_1,\lambda_2,\lambda_ 3):\arctan\lambda_1+\arctan\lambda_2+
\arctan\lambda_3\ge-\arctan c \}, $$
$$\tilde F^3_{c,\Lambda}=\{(\lambda_1,\lambda_2,\lambda_ 3):\arctan\lambda_1+\arctan\lambda_2+
\arctan\lambda_3\ge \pi-\arctan c \}. $$
A simple calculation gives

\medskip
{\bf Lemma 3.2.} {\em $F^i_c$ is strictly, but not uniformly,
elliptic.}

\medskip

Indeed, the  derivatives ${1\over \lambda_i^2+1}>0 $ tend to 0 at
infinity.

\medskip
{\em Remark 3.1.} If we (artificially) impose the uniform
ellipticity condition, we get a smooth solution. Indeed, if ${\bf
F}_{-c}$   verifies this condition on $u$,  the  derivatives
$${1\over \lambda_1^2+1},\;{1\over \lambda_2^2+1},\;
{1\over \lambda_3^2+1}=
{(\lambda_1\lambda_2-1+c\lambda_1+c\lambda_2)^2\over
(1+\lambda_2^2)(\lambda_1^2+1)(c^2+1)}\in\left [{1\over M},M\right
]
$$
for some ellipticity constant $M$, which implies that $u\in
C^{1,1}$ and thus is smooth by [Y].

We give now the principal technical result of this section which
permits to construct in the next section a singular solution of
SLE.

\medskip
 {\bf Proposition 3.1.}{\em There exists a ball
$B=B(0,\varepsilon)$ centered at the origin s.t.

\smallskip
$1).$ The equation
$$\lambda_1\lambda_2+\lambda_2\lambda_3+\lambda_1\lambda_3=\sigma_2(D^2u)=1
$$  has an analytic   solution $u_0$ in $B$ verifying

$(i)$
$$u_0=-{y^4\over
3}+5y^2z^2-x^4+7x^2z^2-{z^4\over 3}+2y^2z-2zx^2+{y^2\over
2}+{x^2\over 2}+O(r^5)
$$

$(ii)$ $$\lambda_1=1+O(r),\; \;\lambda_2= 1+O(r),\;\;\lambda_3=
 -{x^2\over 2}-{3y^2 \over 2}-z^2+O(r^3).$$

$2).$ The equation
$$\lambda_1\lambda_2+\lambda_2\lambda_3+\lambda_1\lambda_3+
 c(\lambda_1\lambda_2\lambda_3-\lambda_1-\lambda_2-\lambda_3)=
 \sigma_2(D^2u)+c(\det(D^2u)-Tr(D^2u))=1$$
 for $c\neq 0, -1$ has an analytic   solution $u_c$ in $B$ verifying

$(i)$

$$u_c={-z^4\over (c+1)(c^2+2c+2)(c^2+c+1)(c^2+1)} +{2z^2y^2(4c^5+4c^4+8c^3+5c^2+4c+4)
\over(c+1)(c^2+c+1)}$$
$$+{2x^2z^2(4c^2+4c+3)  \over (c+1)(c^2+c+1)
}+{y^4(c^2+1)(3c^4+2c^3+2c^2-4c-4) \over (c+1)(c^2+c+1)}
-{x^4(3c^2+2c+2)\over  (c+1)(c^2+c+1) }$$
$$
-2(c^2+1)zy^2+2zx^2+{(c^2+c+1)y^2\over 2}+{(c+1)x^2\over 2}+O(r^5)
$$

$(ii)$

 $$\lambda_1=c^2+c+1+O(r),\; \;\lambda_2= c+1+O(r),\;\;$$
$$\lambda_3=
 -{{{x^2 (c^2+1)(c^2+2c+2) }+{3y^2 c^2(c^2+1)(c^2+2c+2) }+
 {3z^2}\over 2(c+1)(c^2+c+1)(c^2+1)(c^2+2c+2)}}+O(r^3).$$

 3).The equation $(c=-1)$
$$\lambda_1\lambda_2+\lambda_2\lambda_3+\lambda_1\lambda_3+
 c(\lambda_1\lambda_2\lambda_3-\lambda_1-\lambda_2-\lambda_3)=
 \sigma_2(D^2u)-\det(D^2u)+Tr(D^2u)=1$$
  has an analytic   solution $u_{-1}$ in $B$ verifying

$(i)$
$$u_{-1}=48y^2x^2-12y^2z^2-{119x^4\over 2}+93x^2z^2+{z^4\over 2}+2y^2z-9x^2z
-{y^2\over 6}+x^2+O(r^5)$$

 $(ii)$ $$\lambda_1=2+O(r),\;  \;\lambda_2=
 6y^2+6x^2+{3z^2\over 2}+O(r^3),\;\lambda_3= -{1\over 3}+O(r),$$
 for $r=||(x,y,z)||.$}

\medskip
 {\em Proof.} Let us note that $$v_0:=-{y^4\over
3}+5y^2z^2-x^4+7x^2z^2-{z^4\over 3}+2y^2z-2zx^2+{y^2\over
2}+{x^2\over 2},$$
$$v_{-1}:=48y^2x^2-12y^2z^2-{119x^4\over 2}+93x^2z^2+{z^4\over 2}+2y^2z-9x^2z
-{y^2\over 6}+x^2 ,$$ and $$v_c={-z^4\over
(c+1)(c^2+2c+2)(c^2+c+1)(c^2+1)}
+{2z^2y^2(4c^5+4c^4+8c^3+5c^2+4c+4) \over(c+1)(c^2+c+1)}$$
$$+{2x^2z^2(4c^2+4c+3)  \over (c+1)(c^2+c+1)
}+{y^4(c^2+1)(3c^4+2c^3+2c^2-4c-4) \over (c+1)(c^2+c+1)}
-{x^4(3c^2+2c+2)\over  (c+1)(c^2+c+1) }$$
$$
-2(c^2+1)zy^2+2zx^2+{(c^2+c+1)y^2\over 2}+{(c+1)x^2\over 2}
$$ verify their
 respective equations up to the second order, i.e.
 $$\sigma_2(D^2v_0)-1=O(r^3),$$
 $$\sigma_2(D^2v_{-1})-\det(D^2v_{-1})+Tr(D^2v_{-1})-1=O(r^3),$$
 $$\sigma_2(D^2v_c)+c(\det(D^2u_c)-Tr(D^2u_c))-1=O(r^3),$$
which can be proven by a brute force (e. g., MAPLE) calculation
(e.g.,
$$\sigma_2(D^2v_0)-1=4(-10y^4-32y^2x^2-50y^2z^2-42x^4-
130x^2z^2+11z^4-36y^2z+4x^2z+4z^3)$$).

To prove 1). one considers the following Cauchy problem for the
equation $F_0=\sigma_2(D^2u)-1=0\; :$
$${{u}|}_{z=0}={v_0|}_{z=0}= -{y^4\over 3}-x^4+{y^2\over2}+{x^2\over 2},$$
$${\left({\partial u\over \partial z}\right)}_{z=0}=
{\left({\partial v_0\over \partial z}\right)}_{z=0}=
2y^2-2x^2\;.$$

Since the equation is elliptic, we get by the Cauchy-Kowalevskaya
theorem    a unique local analytic solution $u_0$ which should
coincide with $v_0$ within to 4th order.

The same argument is valid for
$$F_c=\sigma_2(D^2u)+c(\det(D^2u)-Tr(D^2u))-1=0$$
 and  the Cauchy problem
$${{u}|}_{z=0}={v_c|}_{z=0}=$$ $$ {y^4(c^2+1)(3c^4+2c^3+2c^2-4c-4)\over 3(c+1)
(c^2+c+1)} -{x^4(3c^2+2c+2)\over 3(c+1)(c^2+c+1) }
+{y^2((c^2+c+1)\over 2}+{x^2(c+1)\over 2} ,$$
$${\left({\partial u\over \partial z}\right)}_{z=0}=
{\left({\partial v_c\over \partial z}\right)}_{z=0}=
-2(c^2+1)y^2+2x^2\;.$$ The claim on the eigenvalues follows
directly from the formulas

$$\det(D^2v_0)= -{x^2\over 2}-{3y^2 \over 2}-z^2+O(r^3),$$
$$\det(D^2v_c)= -{x^2\over 2}-{3y^2 c^2 \over 2}-{3z^2\over(c^2+1)(c^2+2c+2)}+O(r^3)$$
which are straightforward (e.g. $$\det(D^2v_0)=
120y^4x^2-140y^4z^2+168y^2x^4+1440y^2x^2z^2-994y^2z^4-420x^4z^2-1350x^2z^4-
140z^6$$ $$+40y^4z+192y^2x^2z-136y^2z^3-168x^4z-120x^2z^3-16z^5
-10y^4 +20y^2x^2+28y^2z^2-42x^4$$ $$+28x^2z^2-8z^4
-24y^2z+40x^2z-{3y^2\over 2}-{x^2\over 2}-z^2 ).$$ The argument
works for 3). as well.
\section{Legendre Transform}
Let us recall principal properties [see, e.g. CH,  $\S 1.6$] of
the Legendre Transform (for simplicity of notation we consider
here only the case of 3 dimensions used below). Let $f$ be  a
$C^2$ function defined in a domain $D\subset {\bf R}^3$ s.t. its
gradient map $\nabla f: D\longrightarrow {\bf R}^3$ maps
bijectively $D$ onto a domain $G.$ Let $g=(\nabla
f)^{-1}=(P,Q,R):G\longrightarrow D$ be the map inverse to the
gradient. Then the Legendre Transform $ \tilde f:G\longrightarrow
{\bf R}$ is given by
$$  \tilde f(u,v,w):=uP(u,v,w)+vQ(u,v,w)+wR(u,v,w)-f(g(u,v,w)). $$
 Suppose also that $\det(D^2 f)\neq 0$ except for a point $ a\in
 D$ with $b=(\nabla f)(a).$
 Then $(D^2\tilde f)=(D^2 f)^{-1}$ on $G-\{b\}$.

 We want then to apply the Legendre Transform to the solutions
 $u_c$ on a small ball centered at zero.  We need thus to
 verify that   $ \nabla u_c$ is injective. One finds
$$ \nabla u_c=[ U(x,y,z)x+(c+1)x,V(x,y,z)y+(c^2+c+1)y,
 -{4z^3m_c}+x^2W_1(z)+y^2W_2(z)],$$
where $ U(x,y,z),V(x,y,z)\in {\bf R}\{\{x,y,z\}\},
\;W_1(z),W_2(z)\in {\bf R}\{\{z\}\},\;U(0,0,0)=V(0,0,0)=0,\;$
 $\; m_c:={1/ \left((c+1)(c^2+2c+2)(c^2+c+1)(c^2+1)\right)}>0$ for\\ $c\neq
 0,-1$, $m_{-1}:=1/2,$ $m_0:=1/3.$

 Thus the gradient map is injective by [EL, Theorem 1.1] since one has

 \medskip
 {\bf  Lemma 4.1.}
  {\em The ring ${\bf R}\{\{x,y,z\}\}/
  (\partial u_c/\partial x,\partial u_c/\partial y,\partial u_c/\partial z)$ is
   isomorphic to ${\bf R}[h]/(h^3)$.}

\smallskip
{\em Proof.} Indeed, $p:=\partial u_c/\partial x$, $q:=\partial
u_c/\partial y$ can be chosen as new local coordinates and the
lemma follows.

\smallskip
   We can now prove our main result.

 \medskip
 {\bf Theorem 4.1.} {\em Let $\theta \in ]-{\pi\over 2},{\pi\over 2}[$, and let
 $$ {\bf F}_{\theta}(u)=\arctan \lambda_1+\arctan \lambda_2+\arctan \lambda_3-\theta=0.$$
  Then for some ball $ B_{\varepsilon}$
  centered at the origin there exists an analytic function $f_{\theta}$ on
  $\partial  B_{\varepsilon}$ s.t. the unique (Harvey-Lawson) solution $u_{\theta}$
  of the Dirichlet problem  $$\cases{{\bf F}_{\theta}(u)=0 &in $B_{\varepsilon}$\cr
u=f_{\theta} &on $\partial  B_{\varepsilon}\;$\cr}$$ verifies:

\smallskip
$(i)\;\;u_{\theta}\in C^{1,1/3}; $

\smallskip
$(ii)\;\;u_{\theta}\notin C^{1,\delta} $ for $\forall\delta>1/3$.
 }
\medskip

{\em Proof.} We can apply the Legendre Transform to $u_c$ with
$c=\cot(\theta)$ for $\theta\neq 0$, and to $u_0$ for $\theta= 0$
thanks to the injectivity of $ \nabla u_c\;$. Since $u_c$ with
$c\neq 0$ verifies the equation
$$\sigma_2(D^2u)+c(\det(D^2u)-Tr(D^2u))-1=0,$$ its
Legendre Transform  $\tilde u_c$ verifies
$$c(\sigma_2(D^2u)-1)+\det(D^2u)-Tr(D^2u) =0,$$
the signature of $(\lambda_1(\tilde u_c),\lambda_2(\tilde
u_c),\lambda_3(\tilde u_c))$ being (+,+,-) for $c\ge -1$ and
(-,-,+) for $c<-1$ which implies that $\tilde u_c$ lies on the
middle branch of this equation. The same is true for $\tilde u_0$
and the equations
$$\sigma_2(D^2u) -1=0,\;\; \det(D^2u)-Tr(D^2u) =0. $$

The function $\tilde u_c$ is analytic outside zero and  belongs to
$C^{1,1/3}(B_{\varepsilon})$ which proves (i). We need then to
prove that $\tilde u_c$ is a Harvey-Lawson solution  of the
corresponding   Dirichlet problem.

This is implied by   the following form of the Alexandrov maximum principle [A]:

 \bigskip
{\bf Proposition 4.1.} {\em Let $F$ be  a Dirichlet domain,  and
let $u=v+w$ where $v$ is of $\tilde F$-type. If
$$ u\in C^2(B-\{0\})\bigcap C^1(B) $$ and $D^2u$ is non-negatively defined on
 $  B-\{0\}$ then}
$$\sup_B u\le \sup_{\pt B} u.$$
\bigskip

 (ii): Let $u=0,v=0,w\neq 0.$ Then
$$\lambda_3(u_{\theta})=-2m_c^{-1}w^{-2/3}/3 +o(w^{-2/3})$$ which
contradicts the condition $ u_{\theta}\in C^{1,\delta} $ for $
\delta>1/3$.

\medskip
{\em Remark 4.1.} Let us consider the Special Lagrangian
submanifold $L_{u,c} \subset {\bf C}^3$ corresponding to our
singular solution $u_c$ i.e. the graph of the map
$$i\nabla u_c: B\longrightarrow i{\bf R}^3.$$
It is easy to show that it is smooth,  and the singularity of
$u_c$ implies only that the projection $L_{u,c}\longrightarrow B$
is singular (map between smooth manifolds).

\bigskip\bigskip

 \centerline{\bf REFERENCES}

\bigskip\bigskip

\noindent [A]  A.D. Alexandrov; { \it
Some theorems on partial differential equations of the second order}.
Vestnik Leningrad. Univ. 9 (1954), no. 8, 3Ð17.

\medskip
\noindent [Ba] J. Ball; {\it Convexity conditions and existence
theorems in nonlinear elasticity}, Arch. Rat. Mech. Anal. 63
(1977), 337--403.

\medskip
 \noindent [CC] L. Caffarelli, X. Cabre; {\it Fully Nonlinear Elliptic
Equations}, Amer. Math. Soc., Providence, R.I., 1995.

\medskip
 \noindent [CNS] L. Caffarelli, L. Nirenberg, J. Spruck; {\it The Dirichlet
 problem for nonlinear second order elliptic equations III. Functions
  of the eigenvalues of the Hessian, } Acta Math.
   155 (1985), no. 3-4, 261--301.

\medskip
 \noindent [CH]
R. Courant, D. Hilbert; {\it Methods of Mathematical Physics. Vol.
2, Partial Differential Equations},
 Wiley, 1989.

\medskip
 \noindent [CIL]  M.G. Crandall, H. Ishii, P-L. Lions; {\it User's
guide to viscosity solutions of second order partial differential
equations}, Bull. Amer. Math. Soc. (N.S.) 27(1) (1992), 1-67.

\medskip
 \noindent [EL] D. Eisenbud, H. Levin, {\it An algebraic formula for the
  degree of  a $C^{\infty} $ map germ}, Ann. Math. 106 (1977),
  19--44.

\medskip
 \noindent [GT] D. Gilbarg, N. Trudinger; {\it Elliptic Partial
Differential Equations of Second Order, 2nd ed.}, Springer-Verlag,
Berlin-Heidelberg-New York-Tokyo, 1983.

\medskip
\noindent [HL1] R. Harvey, H. B. Lawson Jr., {\it Calibrated
geometries,} Acta Math. 148 (1982), 47--157.

\medskip
 \noindent [HL2] F. R. Harvey, H. B. Lawson Jr.; {\it Dirichlet duality and the
 nonlinear Dirichlet problem} , Comm. Pure Appl. Math. 62(2009), no 3, 396-443.

\medskip
\noindent [NV1] N. Nadirashvili, S. Vl\u adu\c t; {\it
Nonclassical solutions of fully nonlinear elliptic equations,}
Geom. Func. An. 17 (2007), 1283--1296.

\medskip
\noindent [NV2] N. Nadirashvili, S. Vl\u adu\c t; {\it Singular
solutions to fully nonlinear elliptic equations,} J. Math. Pures
Appl. 89 (2008), 107--113.

\medskip
\noindent [NV3] N. Nadirashvili, S. Vl\u adu\c t; {\it On Hessian
fully nonlinear elliptic equations}, arXiv:0805.2694 [math.AP],
submitted.

\medskip
\noindent [T1] N. S. Trudinger; {\it H\"{o}lder gradient estimates
for fully nonlinear elliptic equations,}
  Proc. Roy. Soc. Edinburgh Sect. A  108  (1988),  no. 1-2, 57--65.

\medskip
\noindent [T2] N. S. Trudinger; {\it On regularity and existence
of viscosity solutions of nonlinear second order,
 elliptic equations,} -in :{\it Progr. Nonlinear Differential Equations Appl}.,
 2, Birkh\"{a}user Boston, Boston, MA, 1989,  939--957.

 \medskip
\noindent [T3] N. S. Trudinger; {\it The Dirichlet
 problem for the prescribed curvature equations,}
 Arch. Rat. Mech. Anal. 111 (1990), 153--170.

 \medskip
\noindent [Y] Y. Yuan; {\it A priori estimates for solutions of
fully nonlinear special lagrangian equations,} Ann. Inst. Henri
Pioncar\'{e} non lin\'{e}aire 18 (2001), 261-270.
\end{document}